\documentclass{article}[12pt]
\usepackage{epsfig,makeidx,amsmath,amsfonts,latexsym}
\newtheorem{theorem}{Theorem}[section]
\newtheorem{prop}{Proposition}[section]
\newtheorem{lemma}{Lemma}[section]

\begin{document}

\title{\large\textbf{{Generating stationary random
graphs on $\mathbb{Z}$ with prescribed i.i.d.\ degrees}}}

\author{Maria Deijfen \thanks{Chalmers University of Technology.
E-mail: mia@math.su.se} \and Ronald Meester
\thanks{Vrije Universiteit Amsterdam. E-mail: rmeester@cs.vu.nl }}

\date{February 2006}

\maketitle

\thispagestyle{empty}

\begin{abstract}

\noindent Let $F$ be a probability distribution with support on
the non-negative integers. Two algorithms are described for
generating a stationary random graph, with vertex set
$\mathbb{Z}$, so that the degrees of the vertices are i.i.d.\
random variables with distribution $F$. Focus is on an algorithm
where, initially, a random number of ``stubs'' with distribution
$F$ is attached to each vertex. Each stub is then randomly
assigned a direction, left or right, and the edge configuration is
obtained by pairing stubs pointing to each other, first exhausting
all possible connections between nearest neighbors, then linking
second nearest neighbors, and so on. Under the assumption that $F$
has finite mean, it is shown that this algorithm leads to a
well-defined configuration, but that the expected length of the
shortest edge of a vertex is infinite. It is also shown that any
stationary algorithm for pairing stubs with random, independent
directions gives infinite mean for the total length of the edges
of a given vertex. Connections to the problem of constructing
finitary isomorphisms between Bernoulli shifts are discussed.

\vspace{0.5cm}


\vspace{0.5cm}

\end{abstract}

\section{Introduction}

Recently there has been a lot of interest in the use of random
graphs as models for various types of complex networks. Several
models have been formulated, aiming to capture essential features
of the networks in question such as degree distribution, diameter
and clustering; see for instance Dorogovtsev and Mendes (2003) and
Bollob\'{a}s and Riordan (2003) for surveys. As for the vertex
degree, power-law distributions have been identified in many of
the real-world applications, implying that the ordinary
Erd\H{o}s-Renyi graph, introduced in Erd\H{o}s and Renyi (1959)
and giving Poisson distributed degrees in the limit of large graph
size, is not suitable as a model. This has given rise to a number
of algorithms for generating graphs with an arbitrary prescribed
degree distribution. The most studied one is the so called
configuration model, where each vertex is assigned a random number
of stubs which are then joined pairwise completely at random to
form edges. The asymptotic behavior of this model has been studied
by Molly and Reed (1995,1998), Newman et al.\ (2001) and van der
Hofstad et al.\ (2005) among others. Also, Britton et al.\ (2005)
treats a modification of the model where multiple edges and
self-loops are forbidden, giving simple graphs as a final result.
A different model for generating random graph with prescribed
degrees is studied by Chung and Lu (2002:1,2).

A natural generalization of the problem of generating random
graphs with prescribed degree distributions, is to consider
spatial versions of the same problem, where geometric aspects play
a role. More precisely, given a probability distribution $F$ with
support on the non-negative integers and a set of vertices with
some kind of spatial structure, how should an edge configuration
on this vertex set with degree distribution $F$ be generated?
Clearly the answer of this question depends on the nature of the
spatial structure and also on desired properties of the resulting
configuration.

In this paper we consider the problem of generating a {\em
stationary} random graph, with vertex set $\mathbb{Z}$, and with
degrees that are i.i.d.\ with distribution $F$. We recall that a
graph on $\mathbb{Z}$ is said to be stationary if the distribution
of the edge configuration restricted to any finite subset of
$\mathbb{Z}$ is translation invariant. Hence, we are looking for a
stationary algorithm to obtain edges among the vertices of
$\mathbb{Z}$ in such a way that the vertex degrees become i.i.d.\
random variables with distribution $F$. We have two suggestions of
how to do this.

\subsection{Stepwise pairing with random directions}

Our first suggestion is the one we will spend most of the paper
on. The algorithm runs as follows.

\begin{itemize}
\item[1.] Attach independently a random number of stubs to each
vertex according to the distribution $F$.

\item[2.] Randomly assign a direction, left or right, to each
stub, turning it into an arrow.

\item[3.] Join arrows pointing to each other stepwise, first
exhausting all possible connections between nearest neighbor
vertices, then looking at second nearest neighbors, and so on,
until all arrows are connected.
\end{itemize}

This model will be referred to as the Stepwise Pairing algorithm
with Random Directions (SPRD). In Section 2, we show that, with
probability 1, the SPRD-algorithm leads to a well-defined
configuration, that is, the number of steps required for a right-
(left-) arrow to find a left- (right-) arrow to connect to is
almost surely finite. However, already the shortest edge of a
vertex turns out to have infinite mean; see Section 4. Basically
this follows from properties of a random walk structure that
arises in the analysis of the model. This analysis is complicated
by the fact that the increments of the random walk are not
independent, making standard results inapplicable. In Section 5 we
prove that, if we insist on the directions of the edges of a given
vertex being completely random and also independent of the
configuration at all other vertices, then we can not achieve
finite mean for the {\em total} edge length per vertex in a
stationary way. Note however that, by dropping the requirement
that the directions of the edges should be assigned randomly and
independently, it is in some cases possible to design algorithms
that work; see Examples 5.1 and 5.2.

Readers familiar with Bernoulli isomorphisms might have observed
that our pairing rule is close to the pairing rule used by
Meshalkin (1959) in order to construct an isomorphism between two
specific Bernoulli shifts with equal entropy. In fact, there are
connections, which we will explore a bit in Section 3. We also
mention the paper by Holroyd and Peres (2005), which deals with
stationary matching in a slightly different set-up.

\subsection{Annihilating random walk}

There is a totally different way of generating a stationary graph
with the required properties, making use of random walks and which
we mention for completeness. It can be described in three steps as
follows.

\begin{itemize}
\item[1.] Attach independently a random number of stubs to each
vertex according to the distribution $F$.

\item[2.] To each stub, associate a particle at the same position
on $\mathbb{Z}$ and let all particles start a continuous time
random walk on $\mathbb{Z}$, independent of each other.

\item[3.] Whenever two particles - started at different locations
- meet, draw an edge between the corresponding stubs and remove
the particles from the system.
\end{itemize}

It is not hard to see, and well known (see for instance Arratia
(1981)), that this leads to a limiting configuration in which all
stubs are connected. However, we will not be concerned with this
type of pairing in this paper. See Mattera (2003) for other
connections between annihilating random walks and graphs.

\section{Definition of the SPRD-algorithm}

Let us first describe the SPRD-algorithm in more detail. To begin
with, associate independently to each vertex $i\in\mathbb{Z}$ a
random degree $D_i$ with distribution $F$. Think of this as vertex
$i$ having $D_i$ ``stubs'' sticking out of it. Now turn the stubs
into arrows by randomly associating a direction to each stub. More
precisely, with probability $p$ a stub is pointed to the right and
with probability $1-p$ it is pointed to the left. Write $R_i$
$(L_i)$ for the number of right- (left-) arrows of vertex $i$ and
label the arrows $\{r_{i,j}\}_{j=1}^{R_i}$
$\big(\{l_{i,j}\}_{j=1}^{L_i}\big)$. This gives a configuration
where each vertex $i$ has two ordered sets of arrows
$\{r_{i,j}\}_j$ and $\{l_{i,j}\}_j$ associated to it. These arrows
will now be matched pairwise, a pair always consisting of one
right-arrow and one left-arrow, to create edges between the
vertices. The matching is done stepwise as follows.

\begin{itemize}
\item[1.] First consider all pairs of nearest neighbor vertices
$i$ and $i+1$ and create $\min\{R_i,L_{i+1}\}$ edges between
vertex $i$ and $i+1$ by joining the arrows $r_{i,j}$ and
$l_{i+1,j}$ for $j=1,\ldots,\min\{R_i,L_{i+1}\}$.

\item[2.] Next consider all pairs of second nearest neighbor
vertices $i$ and $i+2$. If, after step 1, there is at least one
unconnected right-arrow at vertex $i$ and at least one unconnected
left-arrow at vertex $i+2$, then we create edge(s) between the
vertices $i$ and $i+2$, by performing all possible connections,
always connecting an arrow $r_{i,j}$ ($l_{i+2,j}$) before
$r_{i,j+1}$ ($l_{i+2,j+1}$).
$$
\vdots
$$
\item[$n$.] In step $n$, we consider all pairs of vertices $i$ and
$i+n$ at distance $n$ from each other and connect arrows that
remain after the previous steps, never using an arrow $r_{i,j}$
($l_{i+n,j}$) before $r_{i,j-1}$ ($l_{i+n,j-1}$).
$$
\vdots
$$
\end{itemize}

The above procedure is clearly stationary, but we have yet to show
that it leads to a well-defined graph. To this end, define the
length of an edge in the resulting configuration to be the
distance between its endpoints. In what follows, we will consider
only the vertex at the origin.

Write $N_{j}^{(r)}$ for the length of the edge created by
right-arrow number $j$ at the origin, $r_{0,j}$, and set
$N_j^{(r)}=\infty$ if $r_{0,j}$ is never connected. Also, define
$N_j^{(l)}$, $j\geq 1$, analogously for the left-arrows. Write
$$
N^{(r)}=\max_{1\leq j\leq R_0}\big\{N_j^{(r)}\big\}\quad
\textrm{and}\quad N^{(l)}=\max_{1\leq j\leq
L_0}\big\{N_j^{(l)}\big\},
$$
and define $N=\max\{N^{(r)},N^{(l)}\}$.

We first show that the algorithm does not work for $p\neq 1/2$.
Here, $P_p$ denotes the probability measure associated with the
SPRD-algorithm when a stub is pointed to the right with
probability $p$.

\begin{prop}\label{prop:transience}
If $p\neq 1/2$, then $P_p(N=\infty)>0$.
\end{prop}

\noindent\textbf{Proof:} To see that $P_p(N_1^{(r)}= \infty)>0$
for $p>1/2$, fix $p>1/2$, let $\Delta_i=L_i-R_{i}$ and define
$$
\left\{
        \begin{array}{l}
        S'_1=L_1;\\
        S'_n=\sum_{i=1}^{n-1}\Delta_i+L_n\textrm{ for }n\geq 2.
        \end{array}
            \right.
$$
Clearly, the first right-arrow at the origin gets connected as
soon as $S_n'$ takes on a positive value and hence it suffices to
show that $P_p(S_n'\leq 0\textrm{ for all }n)>0$. This however
follows easily by noting that the law of large numbers implies
that $S_n'\to\infty$ almost surely.\hfill$\Box$\medskip

Having discarded non-symmetric versions of the algorithm, let us
move on to the symmetric case where the prospects of success
should be better. Indeed, the following proposition guarantees
that, for $p=1/2$, no arrow has to wait infinitely long before it
finds something to connect to.

\begin{theorem}\label{prop:recurrence}
We have $P_{1/2}(N<\infty)=1$.
\end{theorem}

\noindent \textbf{Proof:} For ease of notation, write $P_{1/2}=P$.
First note that, by symmetry, it suffices to show that
$P(N^{(r)}<\infty)=1$. Also, by the definition of the algorithm,
the arrows $\{r_{0,j}\}$ are used in chronological order, implying
that $N_j^{(r)}\leq N_{j+1}^{(r)}$. It follows that
$N^{(r)}=N_{R_0}^{(r)}$, and, since $R_0<\infty$ almost surely, we
are done if we can show that $P(N_j^{(r)}<\infty)=1$ for all
$j$.

To do this, we first consider the case $j=1$ and show that
$P(N_1^{(r)}<\infty)=1$, that is, the length of the edge created
by the first right-arrow $r_{0,1}$ at the origin is almost surely
finite. This is done by dominating the length of the edge by the
time at which a recurrent random walk takes on a positive value
for the first time. To be more specific, define $\Delta_i=L_i-R_i$
and write $S_n=\sum_{i=1}^n\Delta_i$. The variables $\{\Delta_i\}$
are i.i.d.\ and symmetric, implying that $\eta:=\inf\{n;S_n>0\}$
is finite with probability 1. Now note that, as soon as $S_n>0$,
we know that the arrow $r_{0,1}$ must have found a left-arrow to
connect to. Indeed, if $S_n>0$, we also have $S_n+R_n>0$, and the
fact that $S_n+R_n>0$ means that the total number of left-arrows
on the vertices $1,\ldots,n$ is strictly larger than the total
number of right-arrows on the vertices $1,\ldots,n-1$, implying
that, at some vertex $1,\ldots,n$, there must be a left-arrow for
$r_{0,1}$ to connect to. It follows that $N_1^{(r)}\leq \eta$ and
we are done.

Now assume in an inductive fashion that $P(N_j^{(r)}<\infty)=1$
and suppose that $P(N_{j+1}^{(r)}=\infty)>0$. Write
$\Psi=\{\Psi_i\}=\{(L_i,R_i)\}$ for the random configuration of
arrows at the vertices and, for configurations with
$N_j^{(r)}<\infty$, introduce a coupled configuration
$\widehat{\Psi}=\{\widehat{\Psi}_i\}$ that is identical to $\Psi$
except that the directions of the stubs at the vertex $N_j^{(r)}$
are generated independently. Let $\widehat{N}_j^{(r)}$ be the
length of the edge formed by $r_{0,j}$ in $\widehat{\Psi}$ and
define
$$
A_j=\left\{N_{j+1}^{(r)}=\infty\right\}\cap\left\{\widehat{L}_{N_j^{(r)}}=0\right\}.
$$
Note that, on the event $A_j$, we have
$\widehat{N}_j^{(r)}=\infty$: Indeed $r_{0,j}$ cannot connect
before vertex $N_j^{(r)}$ in $\widehat{\Psi}$, since left-arrows
have been removed at $N_j^{(r)}$ without any new right-arrows
being added. Furthermore, if $r_{0,j}$ is connected to a
left-arrow at vertex $m\geq N_j^{(r)}$ in $\widehat{\Psi}$ it
would imply that $r_{0,j+1}$ was connected at the latest to $m$ in
$\Psi$ and this conflicts with the fact that
$N_{j+1}^{(r)}=\infty$. Hence we have
$P(\widehat{N}_{j}^{(r)}=\infty)\geq P(A_j)$. It follows from
the assumption that $P(A_j)>0$ and, since clearly $N_j^{(r)}$ and
$\widehat{N}_j^{(r)}$ have the same distribution, we have shown
that $P(N_{j}^{(r)}=\infty)>0$. But this is a contradiction and,
by induction over $j$, we conclude that
$P(N_j^{(r)}<\infty)=1$ for all $j$, as desired.\hfill$\Box$

\section{Connections to Bernoulli isomorphisms}

Consider a stochastic process $X$ indexed by ${\mathbb Z}$ with
i.i.d.\ marginals taking values in $1,2,\ldots,s$, with
probabilities $p_1,\ldots, p_s$ respectively. The process $X$ is
often called a Bernoulli shift, and is identified with the vector
$(p_1,\ldots, p_s)$. Next, consider another such process $Y$, with
values in $1,2,\ldots, t$ and probabilities $q_1,\ldots, q_t$,
respectively. We write $S_X=\{1,2,\ldots, s\}^{\mathbb Z}$ and
$S_Y=\{1,2,\ldots, t\}^{\mathbb Z}$. Loosely speaking, $X$ and $Y$
are called {\em ismorphic} if there exists a pairing of almost all
realisations of $X$ and $Y$ in a bijective way, such that the
pairing commutes with the shift operator on $S_X$ and $S_Y$.

Meshalkin (1959) was one of the first to explicitly identify such
a coding between two particular Bernoulli shifts, namely between
$(\frac14, \frac14, \frac14, \frac14)$ and $(\frac12, \frac18,
\frac18,\frac18, \frac18)$. His coding corresponds to our
algorithm in the case when each vertex has at most one edge
associated to it, as follows. Associate to each edge a random
label, $a$ or $b$, independently and with equal probability. This
leads to four equally likely symbols, namely $(l,a), (l,b), (r,a)$
and $(r,b)$, where $l$ and $r$ refer to the edge pointing to the
left or to the right. The coding is now defined so that, whenever
we see $(r,a)$ or $(r,b)$ we write an $r$, and whenever we see an
$l$ we write $(l,x,y)$, where $x$ and $y$ are the symbols
corresponding to the edge that is formed with the unique stub at
that position. It is not hard to see that this codes the original
four symbols into five new symbols with probabilities $\frac12$
and four times $\frac18$, and that this coding is invertible.
Furthermore, the coding is {\em finitary}, that is, one has to
look only a (random) finite distance in both directions to see
what symbol that should be written in the coding. Indeed, once we
have identified the stub that connects to our current stub of
interest, we can write down the correct symbol.

This idea can be stretched to apply for a general degree
distribution $F$ with bounded support in our pairing algorithm,
that is, every $F$ with bounded support leads to an isomorphism
between two particular Bernoulli shifts, as the reader can easily
verify. Certain results of coding between Bernoulli shifts then
have corollaries for our algorithm. We mention the well known fact
(see for instance Parry (1979) or Schmidt (1984)) that in any
nontrivial situation, the expected distance that has to be
explored in a finitary coding between two Bernoulli shifts with
equal entropy has infinite expectation. From this it follows that
the {\em longest} edge at a given vertex in the SPDR-algorithm has
infinite expected length. Below we strengthen this result to the
shortest edge, which does not have an interpretation in the coding
setup described here.

\section{The mean length of the shortest edge}

For the remainder of the paper we only consider the symmetric
SPRD-algorithm, which, by Theorem \ref{prop:recurrence}, leads to
well-defined configurations. The next task is to look at the
expected length of the edges. For distributions with finite mean,
we will prove the following theorem.

\begin{theorem}\label{prop:mean}
If $F$ has finite mean, then both ${\rm{E}}\big[N^{(r)}_1\big]$
and ${\rm{E}}\big[N^{(l)}_1\big]$ are infinite.
\end{theorem}

Define $X_i=L_i-R_{i-1}$, that is, $X_i$ is a difference between
the number of left and right arrows involving two
\emph{neighboring} vertices. Let
$S_n^{(m)}=\sum_{i=m+1}^{m+n}X_i$, and write
$\tau_{\uparrow}^{(m,x)}$ for the first time when the process
$S_n^{(m)}$ reaches above the level $x$, that is,
$$
\tau_{\uparrow}^{(m,x)}=\min\{n;S_n^{(m)}\geq x\}.
$$
Clearly, to prove Theorem \ref{prop:mean}, it suffices to show
that E$[N^{(r)}_1]=\infty$. To see why this should be the case,
note that, if $L_1=0$ (which happens with positive probability),
then the first right-arrow at the origin is connected when
$S_n^{(1)}$ takes on a value larger than or equal to 0, that is,
$N^{(r)}_1=\tau_{\uparrow}^{(1,0)}+1$. If $S_n^{(1)}$ had
independent increments, it would follow from standard random walk
theory that $\tau_{\uparrow}^{(1,0)}$ had infinite mean. However,
$X_i$ and $X_{i+1}$ are not independent, since information about
the arrow configuration at vertex $i$ is used for both variables.

Let $\mu$ denote the mean of $F$. The following lemma will play a
key role in the proof of Theorem \ref{prop:mean}.

\begin{lemma}\label{lemma}
For all $i\in\mathbb{Z}$, we have
{\rm{E}$[\tau_{\uparrow}^{(i,2\mu)}]=\infty$}.
\end{lemma}

\noindent \textbf{Proof of Lemma \ref{lemma}:} By stationarity, it
suffices to show that E$[\tau_{\uparrow}^{(0,2\mu)}]=\infty$.
Assume for contradiction that
$\textrm{E}[\tau_{\uparrow}^{(0,2\mu)}]<\infty$ and define
$\tau_{\downarrow}^{(m,x)}$ to be the first time when the process
$S_n^{(m)}$ reaches below the level $x$, that is,
$$
\tau_{\downarrow}^{(m,x)}=\min\{n;S_n^{(m)}\leq x\}.
$$
Note that, by symmetry, we have
$\textrm{E}[\tau_{\downarrow}^{(0,-2\mu)}]=\textrm{E}[\tau_{\uparrow}^{(0,2\mu)}]$.
The idea of the proof is to use the finite mean assumption to
create a linear negative drift for the process $S_n^{(0)}$. By
symmetry, $S_n^{(0)}$ must then also have the same positive drift
and to maintain both these drifts it is forced to oscillate more
and more vigorously between large positive and large negative
values, something which it will not be able to do in the long run.
To turn this heuristics into a proof, introduce an i.i.d.\
sequence $\{\Delta\tau_j\}$ with mean
$\textrm{E}[\tau_{\downarrow}^{(0,-2\mu)}]+1$ by defining
$$
\left\{\begin{array}{l}
       \Delta\tau_0=0;\\
       \Delta\tau_j=\min\left\{n;S_n^{\left(\sum_{i=0}^{j-1}\Delta\tau_i\right)}
       \leq -2\mu\right\}+1,\textrm{ for }j\geq 1;
       \end{array}
\right.
$$
here, the +1 is added to get independence, noting that $X_i$ and
$X_k$ are independent as soon as $|i-k|\geq 2$. This sequences
gives rise to a renewal process with time increments
$\{\Delta\tau_j\}$ and events referred to as
\emph{down-transitions} occurring at the time points $\{\tau_i\}$,
where $\tau_i=\sum_{j=1}^i\Delta\tau_j$. Write $M_n$ for the
number of down-transitions in the time interval $[0,n]$ and note
that, by the renewal theorem, we have
$$
\frac{M_n}{n}\longrightarrow\big(\textrm{E}[\tau^{(0,-2\mu)}]+1\big)^{-1}\quad\textrm{
a.s.
as }n\rightarrow\infty.
$$
Hence, defining $2c=\big(\textrm{E}[\tau^{(0,-2\mu)}]+1\big)^{-1}$
and $E_m=\{M_n>nc\textrm{ for all }n\geq m\}$, it follows that

\begin{equation}\label{eq:E_m}
P(E_m)\rightarrow 1\quad\textrm{as }m\rightarrow\infty.
\end{equation}

At the point $\tau_d$ of the $d$-th down-transition we have
$$
S_{\tau_d}^{(0)}\leq -2\mu d+\sum_{i=1}^d X_{\tau_i},
$$
where $X_{\tau_i}=L_{\tau_i}-R_{\tau_i-1}\leq L_{\tau_i}$. The
degree of a vertex $\tau_i-1$ is atypical, since it is defined as
a first passage point for the process $S_n^{(\tau_{i-1})}$.
However, the vertex $\tau_i$ has the unconditional degree
distribution $F$, meaning that E$[L_{\tau_i}]=\mu/2$. Also, since
$|\tau_i-\tau_{i-1}|\geq 2$, the variables $\{L_{\tau_i}\}$ are
independent. Combining this we get from the strong law of large
numbers that
$$
\frac{1}{d}\sum_{i=1}^dX_{\tau_i}\leq
\frac{1}{d}\sum_{i=1}^nL_{\tau_i}\rightarrow\frac{\mu}{2}\quad\textrm{a.s.
as }n\rightarrow\infty,
$$
and, defining
$$
F_m=\bigg\{\sum_{i=1}^{\lfloor nc\rfloor}X_{\tau_i}\leq \lfloor
nc\rfloor\mu\textrm{ for all }n\geq m\bigg\},
$$
it follows that

\begin{equation}\label{eq:F_m}
P(F_m)\rightarrow 1\textrm{ as }m\rightarrow \infty.
\end{equation}

Note at this point that, if the sequence ${\tau_i}$ were to be
defined in terms of up-transitions instead of down-transitions,
then, to estimate the value of the process after some large number
of transitions would require a \emph{lower} bound for the sum of
the auxiliary steps $X_{\tau_i}$. This however would cause
trouble, since, as mentioned above, the negative part of a
variable $X_{\tau_i}$ concerns the arrow configuration at a first
passage vertex which is presumably difficult to control. Hence we
are in the peculiar situation of being able to show a statement
for down-transitions but not for up-transitions directly, in an
otherwise completely symmetric situation.

Next divide $\mathbb{Z}^+$ into intervals
$\{\mathcal{I}_k\}_{k\geq 0}$ of length $l$, where
$\mathcal{I}_k=\{i;kl\leq i<(k+1)l\}$, and write $B_k$ for the
event that the interval $\mathcal{I}_k$ contains a down-transition
in the renewal process $\{\tau_i\}$. Clearly, by picking $l$
large, we can make sure that $P(B_k)\geq 0.99$ for all $k$. Define
$Y_k$ to be the sum of the degrees of all vertices in
$\mathcal{I}_k$, that is, $Y_k=\sum_{i=kl}^{(k+1)l-1}D_i$. The
distribution of $Y_k$ does not depend on $k$ and $Y_k<\infty$
almost surely, implying that

\begin{equation}\label{eq:Y_k}
P(Y_k\geq 2\mu\lfloor klc\rfloor)\rightarrow 0 \textrm{ as
}k\rightarrow\infty.
\end{equation}

By (\ref{eq:E_m}), (\ref{eq:F_m}) and (\ref{eq:Y_k}), if we pick
$k$ large enough, we have

\begin{itemize}
\item[(i)] $P(E_{kl})\geq 0.99$; \item[(ii)] $P(F_{kl})\geq
0.99$; \item[(iii)] $P(Y_k\geq 2\mu\lfloor klc\rfloor)\leq 0.5$.
\end{itemize}

Fix such a $k$ and define
$$
D_k^-=\{\exists n\in\mathcal{I}_k \textrm{ such that
}S_n^{(0)}\leq -\mu\lfloor klc\rfloor\}
$$
and
$$
D_k^+=\{\exists n\in\mathcal{I}_k \textrm{ such that
}S_n^{(0)}\geq \mu\lfloor klc\rfloor\}.
$$
Now observe that $B_k\cap E_{kl}\cap F_{kl}\subset D_k^-$: Indeed,
the event $E_{kl}$ implies that $m\geq \lfloor klc\rfloor$
down-transitions have occurred in $[0,kl]$, and $F_{kl}$ implies
that $\sum_{i=1}^{m+1}X_{\tau_i}\leq \mu (m+1)$. Hence, at the
point $\tau_{m+1}$ of the next down-transition, we have

\begin{eqnarray*}
S_{\tau_{m+1}}^{(0)} & \leq & -2\mu (m+1)+\sum_{i=1}^{m+1}X_{\tau_i}\\
& \leq & -2\mu\lfloor klc\rfloor+\mu\lfloor klc\rfloor\\
& = & -\mu\lfloor klc\rfloor.
\end{eqnarray*}

\noindent But this means that $D_k^-$ must occur, since, on the
event $B_k$, at least one down-transition is to take place in
$\mathcal{I}_k$, that is, $\tau_{m+1}\in \mathcal{I}_k$. It
follows that

\begin{eqnarray*}
P(D_k^-) & \geq & P(B_k\cap E_{kl}\cap F_{kl})\\
& \geq & 1-P(B_k^c)-P(E_{kl}^c)-P(F_{kl}^c)\\
& \geq & 0.97.
\end{eqnarray*}

By symmetry, we have $P(D_k^+)=P(D_k^-)$ and hence $P(D_k^+\cap
D_k^-)\geq 0.94$. Now note that, on $D_k^+\cap D_k^-$, we are to
visit both a state above the level $\mu\lfloor klc\rfloor$ and a
state below the level $-\mu\lfloor klc\rfloor$ in the interval
$\mathcal{I}_k$, meaning that $Y_k\geq 2\mu\lfloor klc\rfloor$ on
$D_k^+\cap D_k^-$. Thus

\begin{eqnarray*}
P(Y_k\geq 2\mu\lfloor klc\rfloor) & \geq & P(D_k^+\cap D_k^-)\\
& \geq & 0.94.
\end{eqnarray*}

\noindent But this contradicts (iii) in the choice of $k$. Hence
the assumption that
$\textrm{E}[\tau_{\uparrow}^{(0,2\mu)}]<\infty$ must fail and the
lemma is proved.\hfill$\Box$\medskip

\noindent \textbf{Proof of Theorem \ref{prop:mean}:} By symmetry,
it suffices to show that E$[N_1^{(r)}]=\infty$. To do this, as
before, write $\Psi=\{\Psi_i\}=\{(L_i,R_i)\}$ for the random
configuration of arrows at the vertices and pick $k$ so large that
$P(\sum_{i=1}^kD_i\geq 2\mu)>0$. Introduce a coupled configuration
$\widehat{\Psi}=\{\widehat{\Psi}_i\}$ with the same degrees at all
vertices and $\widehat{\Psi}_i=\Psi_i$ for
$i\not\in\{1,\ldots,k+1\}$, but where the directions of the arrows
at the vertices $1,\ldots,k+1$ are generated independently. Define
$$
A=\left\{\sum_{i=1}^kD_i\geq
2\mu\right\}\cap\left\{\widehat{L}_i=0\textrm{ for all
}i=1,\ldots,k+1\right\}
$$
and let $\widehat{N}_{1}^{(r)}$ be the length of the edge formed
by $r_{0,1}$ in $\widehat{\Psi}$. We then have
$$
\textrm{E}\big[\widehat{N}_{1}^{(r)}\big]\geq
\textrm{E}\big[\widehat{N}_{1}^{(r)}\big|A\big]P(A).
$$
Since clearly $\widehat{N}_{1}^{(r)}$ has the same distribution as
$N_{1}^{(r)}$ and $P(A)>0$, we are done if we can show that
$\textrm{E}\big[\widehat{N}_{1}^{(r)}\big| A\big]=\infty$. To this
end, let $\widehat{S}_n^{(m)}$ be defined in the same way as
$S_n^{(m)}$ but based on the coupled configuration
$\widehat{\Psi}$ and write
$\hat{\tau}_{\uparrow}^{(m,x)}=\inf\{n;\widehat{S}_n^{(m)}\geq
x\}$. On $A$, there are in total at least $2\mu$ right-arrows
attached to the vertices $1,\ldots ,k$ while there are no
left-arrows at all on the vertices $1,\ldots,k+1$. Thus, a
right-arrow at the origin can not be connected until the process
$\widehat{S}_n^{(k+1)}$ takes on a value larger than $2\mu$. It
follows that
$$
\textrm{E}\big[\widehat{N}_{1}^{(r)}\big|A\big]\geq
(k+1)+\textrm{E}[\hat{\tau}_{\uparrow}^{(k+1,2\mu)}|A].
$$
The effect that the conditioning on $A$ has on
$\hat{\tau}_{\uparrow}^{(k+1,2\mu)}$ is that the first term in the
unconditional sum $\widehat{S}_n^{(k+1)}$ is replaced by
$L_{k+2}-D_{k+1}$, since, on $A$, all $D_{k+1}$ stubs at vertex
$k+1$ point to the right. This means that, conditional on $A$, the
passage time $\hat{\tau}_{\uparrow}^{(k+1,2\mu)}$ is
stochastically larger than in the unconditional case, implying
that E$[\hat{\tau}_{\uparrow}^{(k+1,2\mu)}|A]\geq
\textrm{E}[\hat{\tau}_{\uparrow}^{(k+1,2\mu)}]$. Hence
$$
\textrm{E}\big[\widehat{N}_{1}^{(r)}|A\big]\geq
(k+1)+\textrm{E}[\hat{\tau}_{\uparrow}^{(k+1,2\mu)}].
$$
Since $\hat{\tau}_{\uparrow}^{(k+1,2\mu)}$ has the same
distribution as $\tau_{\uparrow}^{(k+1,2\mu)}$, it follows from
Lemma \ref{lemma} that
E$[\hat{\tau}_{\uparrow}^{(k+1,2\mu)}]=\infty$ and the theorem is
proved.\hfill$\Box$

\section{Finite mean is impossible}

We are now at the point of having formulated a stationary
algorithm that takes a discrete distribution $F$ as input and
produces a stationary random edge configuration on $\mathbb{Z}$
with i.i.d.\ vertex degrees with distribution $F$. Provided that
$F$ has finite mean, all connections are almost surely finite but
the expected length of the connections is infinite. The obvious
question is: Can we do better? The following simple examples show
that, if we no longer assign i.i.d.\ directions to the stubs,
then, for certain distributions $F$, indeed we can.

\medskip\noindent \textbf{Example 5.1} Write $f_j$ for the probability that
a given vertex has degree $j$, fix $n\in\mathbb{N}$ and let $F$ be
defined by
$$
\left\{
\begin{array}{ll}
            f_j=(n+1)^{-1} & j=0,2,4,\ldots,2n;\\
            f_j=0 & j\not\in \{0,2,\ldots,2n\}.
                    \end{array}
            \right.
$$
A configuration with this degree distribution and connections with
finite mean is generated by proceeding in the same way as in the
SPRD-algorithm except that the directions of the stubs are not
assigned randomly but according to the deterministic rule that a
vertex with degree $2k$ is equipped with exactly $k$ arrows in
each direction. To see this, note that, assuming that the origin
has degree $d$, all right-arrows at the origin will be connected
as soon as a vertex $i\geq 1$ with degree larger than $d$ is
encountered. The expected distance until we come across a vertex
with degree exactly $d$ is $f_d^{-1}$ and removing the
conditioning on $d$ it follows that the expected length of the
longest connection to the right is bounded by $n$. By symmetry,
the expected maximal length to the left is also bounded by
$n$.\hfill$\Box$\medskip

\noindent \textbf{Example 5.2} Let $F=\delta_1$, that is, every
vertex is to have exactly one edge connected to it. To generate
such a configuration, attach one stub to each vertex and then
imagine that a coin is flipped. If the coin comes up heads the
stubs at the odd vertices are pointed to the right and the stubs
at the even vertices to the left and if it comes up tails we do
the other way around. The arrows are then connected according to
the stepwise paring algorithm. It is easy to see that with this
procedure (which is clearly stationary), all connections will end
up having length 1.\hfill$\Box$\medskip

Recall that, in the SPRD-algorithm, the directions of the stubs
are assigned (i) randomly, and (ii) independently, for each stub.
This gives rise to a random walk type structure which is recurrent
but has infinite mean. In Example 5.1 above, the directions of the
edges are not random and, in Example 5.2, they are not
independent. This destroys the random walk arguments and makes it
possible to obtain configurations where the connections have
finite mean. Thus, for some distributions $F$, it is indeed
possible to outdo the SPRD-algorithm by being clever when
assigning the directions of the stubs. However, we conjecture
that, if the directions are assigned independently for each stub,
then it is impossible to formulate a rule for connecting
right-arrows with left-arrows so that the expected length of the
resulting edges becomes finite. A weaker formulation of this
conjecture is proved in Theorem \ref{prop:greedy} below.

Let $\Psi$ be a random configuration of arrows on $\mathbb{Z}$
generated by the RD-algorithm, that is, first a random number of
stubs with distribution $F$ is attached to each vertex and then
each stub is randomly assigned the direction left or right. An
algorithm $\mathcal{A}$ for connecting the arrows in $\Psi$ will
be called a \emph{pairing rule} if, with probability 1, each
left-arrow is connected to exactly one right-arrow and each
right-arrow is connected to exactly one left-arrow. Furthermore,
$\mathcal{A}$ is said to be stationary if the resulting joint edge
length distributions are translation invariant. For a given
pairing rule $\mathcal{A}$, write $T_{\mathcal{A}}$ and
$N_{\mathcal{A}}$ for the total length of all edges connected to
the origin and the length of the longest edge connected to the
origin respectively.

\begin{theorem}\label{prop:greedy}
If $F$ has finite mean, then, for all stationary pairing rules
$\mathcal{A}$, we have that ${\rm{E}}[T_{\mathcal{A}}]=\infty$.
If, in addition, $F$ has bounded support, then
${\rm{E}}[N_{\mathcal{A}}]=\infty$.
\end{theorem}

The proof of this theorem is based on a combinatorial lemma
involving the concept of \emph{nested} graphs. To define this
concept, consider a given edge configuration
$\{(i,j)\}_{i,j\in\mathbb{Z}}$ on $\mathbb{Z}$. Two edges $(i,j)$
and $(i',j')$ are said to {\em cross} each other if $i<i'<j<j'$ or
$i'<i<j'<j$ and the configuration $\{(i,j)\}_{i,j\in\mathbb{Z}}$
is called {\em nested} if it does not contain any crossing edges.
An important observation is that, for a given configuration $\psi$
of arrows on $\mathbb{Z}$, there is a unique nested edge
configuration, to be denoted by $\mathcal{N}_{\psi}$, which is
obtained by the stepwise pairing algorithm. Indeed, to avoid
crossing edges we are forced to perform all possible connections
between vertices at distance $n=1,2,\ldots$, starting with $n=1$,
and, conversely, successively performing all possible connections
between vertices at distance $n$, with $n$ increasing, can never
in any step create crossing edges, since this would mean that a
possible connection in a previous step was missed.

To formulate the aforementioned lemma, write $\Gamma$ for the set
of all arrow configurations $\psi$ on $\mathbb{Z}$ for which all
edges in $\mathcal{N}_{\psi}$ are finite. Pick $\psi\in\Gamma$
and, for an edge $e\in\mathcal{N}_{\psi}$, let $\psi_e^{(r)}$ and
$\psi_e^{(l)}$ be the set of right-arrows and left-arrows
respectively in $\psi$ that are used to form the edge $e$ and the
edges `under' $e$ in $\mathcal{N}_{\psi}$. More precisely, if $e$
is made up of the arrows $r_{i,j}$ and $l_{i+n,j'}$, then
$\psi_e^{(r)}$ consists of the arrows $\{r_{i,k}\}_{k=1}^j$
together with all right-arrows at the vertices $i+1,\ldots
,i+n-1$, and $\psi_e^{(l)}$ consists of $\{l_{i+n,k}\}_{k=1}^{j'}$
and all left-arrows at the vertices $i+1,\ldots ,i+n-1$. Write
$t_e(\mathcal{N}_{\psi})$ for the total length of all edges
`under' $e$ in ${\cal N}_{\psi}$.

Next, let $\mathcal{E}_\psi$ be some edge configuration
based on the same arrow configuration $\psi$. Call an edge in
$\mathcal{E}_\psi$ a $\psi_e^{(r)}$-\emph{edge} if it contains an
arrow belonging to the set $\psi_e^{(r)}$ and let
$t_e^{(r)}(\mathcal{E}_{\psi})$ denote the total length of all
$\psi_e^{(r)}$-edges in the configuration $\mathcal{E}_{\psi}$.
Define $t_e^{(l)}(\mathcal{E}_{\psi})$ analogously. The lemma now
reads as follows.

\begin{lemma}\label{lemma:nest}
For all $\psi\in\Gamma$, all configurations ${\cal E}_{\psi}$
based on $\psi$, and all $e\in\mathcal{N}_\psi$, we have
$t_e(\mathcal{N}_\psi)\leq t_e^{(r)}(\mathcal{E}_\psi)$ and
$t_e(\mathcal{N}_\psi)\leq t_e^{(l)}(\mathcal{E}_\psi)$.
\end{lemma}

\noindent\textbf{Proof of Lemma \ref{lemma:nest}:} Fix a
$\psi\in\Gamma$, an edge $e\in\mathcal{N}_\psi$ and an edge
configuration $\mathcal{E}_\psi$ based on $\psi$. Define
$w_k^{(r)}$ to be the number of $\psi_e^{(r)}$-edges in
$\mathcal{E}_\psi$ that crosses the interval $[k-1,k]$. More
precisely, $w_k^{(r)}$ is the number of edges in
$\mathcal{E}_\psi$ that has its left endpoint at a vertex $l\leq
k-1$, its right endpoint at $l'\geq k$ and that is created by a
right-arrow that belongs to $\psi_e^{(r)}$. Also, let
$\widetilde{w}_k^{(r)}$ be the same quantity in the nested
configuration $\mathcal{N}_\psi$. We will show that

\begin{equation}\label{eq:m_k}
w_k^{(r)}\geq \widetilde{w}_k^{(r)}\textrm{ for all }k.
\end{equation}

\noindent Since clearly
$t_e(\mathcal{N}_\psi)=\sum_{k=-\infty}^\infty
\widetilde{w}_k^{(r)}$ and
$t_e^{(r)}(\mathcal{E}_\psi)=\sum_{k=-\infty}^\infty w_k^{(r)}$
this implies that $t_e(\mathcal{N}_\psi)\leq
t_e^{(r)}(\mathcal{E}_\psi)$. The inequality
$t_e(\mathcal{N}_\psi)\leq t_e^{(l)}(\mathcal{E}_\psi)$ is proved
similarly.

To establish (\ref{eq:m_k}), assume that the edge $e$ connects the
vertices $i$ and $i+n$, and is created by right-arrow number $j$
at vertex $i$. In the nested configuration, all arrows in
$\psi_e^{(r)}$ are connected to left-arrows at the vertices
$i+1,\ldots,i+n$, meaning that $\widetilde{w}_k^{(r)}=0$ for
$k\not\in\{i+1,\ldots,i+n\}$, and hence trivially $w_k^{(r)}\geq
\widetilde{w}_k^{(r)}$ for such $k$. To deal with
$k\in\{i+1,\ldots,i+n\}$, note that in any edge configuration
based on $\psi$, at least $j$ $\psi_e^{(r)}$-edges must cross the
interval $[i,i+1]$, implying that $w_{i+1}^{(r)}\geq j$.
Furthermore, the interval $[i+1,i+2]$, must be crossed by at least
$j+r_{i+1}-l_{i+1}$ $\psi_e^{(r)}$-edges, where $r_{i+1}$
($l_{i+1}$) denotes the number of right- (left-) arrows at vertex
$i+1$. Hence $w_{i+2}^{(r)}\geq j+r_{i+1}-l_{i+1}$. Continuing in
the same way, we obtain lower bounds for all $w_k^{(r)}$'s with
$k\in\{i+1,\ldots,i+n\}$ and, from the construction of the nested
configuration $\mathcal{N}_\psi$, it follows that these bounds
hold with equality for the $\widetilde{w}_k^{(r)}$'s, and
(\ref{eq:m_k}) follows. \hfill$\Box$\medskip

\noindent \textbf{Proof of Theorem \ref{prop:greedy}:} Let
$\mathcal{A}$ be a stationary pairing rule for an arrow
configuration $\Psi$ generated by the RD-algorithm. If, with
positive probability, $\mathcal{A}$ gives rise to configurations
with infinitely long connections, the conclusion of the
proposition is immediate. Thus assume that all edges in a
configuration obtained from $\mathcal{A}$ are finite almost surely,
and write $T_{\mathcal{A}}^{(r)}$ and $T_{\mathcal{A}}^{(l)}$ for
the total length of the edges created by the right-arrows and
left-arrows respectively at the origin in an edge configuration
generated by $\mathcal{A}$. We will show that E$[T_{\mathcal{A}}^{(r)}]$
and E$[T_{\mathcal{A}}^{(l)}]$ are both infinite.

To prove that E$[T_{\mathcal{A}}^{(r)}]=\infty$, let
$T_{\mathcal{N}}^{(r)}$ be the total length of all edges created
by the right-arrows at the origin in a nested configuration
obtained from the stepwise pairing algorithm and note that, by
Theorem \ref{prop:mean}, we have
E$[T_{\mathcal{N}}^{(r)}]=\infty$. If, with probability 1,
$\mathcal{A}$ results in a nested configuration, then
$T_{\mathcal{A}}^{(r)}$ has the same distribution as
$T_{\mathcal{N}}^{(r)}$ and the claim follows. So assume that with
positive probability $\mathcal{A}$ produces unnested
configurations and let $\mathcal{E}$ be such a configuration with
underlying arrow configuration $\psi$. Write $t_i^{(r)}$ and
$\tilde{t}_i^{(r)}$ for the total length of the edges created by
the right-arrows at vertex $i$ in the configuration $\mathcal{E}$
and $\mathcal{N}_{\psi}$ respectively and let
$\widetilde{m}_i^{(r)}$ be the length of the longest edge formed
by the right-arrows at vertex $i$ in $\mathcal{N}_\psi$. It
follows from Lemma \ref{lemma:nest} that, for all $i$, we have
$$
\sum_{j=i}^{i+\widetilde{m}_i^{(r)}-1}t_j^{(r)}\geq
\sum_{j=i}^{i+\widetilde{m}_i^{(r)}-1}\tilde{t}_j^{(r)}.
$$
By the ergodic theorem, E$[T_{\cal N}^{(r)}]$ is equal to the
average of $\tilde{t}_j^{(r)}$ and hence, for every realization of
${\cal A}$, the average right-degree is bounded below by
E$[T_{\cal N}^{(r)}]$, proving that E$[T_{\cal A}^{(r)}]=\infty$.
That E$[T_{\cal A}^{(l)}]=\infty$ is proved analogously and the
first claim of the theorem follows.

The second claim is established by noting that if $k$ is an upper
bound for the support of $F$, we have $T_{\mathcal{A}}\leq
kM_{\mathcal{A}}$.\hfill$\Box$

\section*{References}

\noindent Arratia, R. (1981), Limiting point processes for rescalings of
coalescing and annihilating random walks on ${\mathbb Z}^d$, {\em The
Annals of Probability} {\bf 9}, 909-936.\medskip

\noindent Bollob\'{a}s, B. and Riordan, O. (2003), Mathematical
results on scale-free random graphs, in \emph{Handbook of graphs
and networks}, Wiley, pp 1-34.\medskip

\noindent Britton, T., Deijfen, M. and Martin-L\"{o}f, A. (2005),
Generating simple random graphs with prescribed degree
distribution, preprint (www.math.su.se/$\sim$mia).\medskip

\noindent Chung, F. and Lu, L. (2002:1), Connected components in
random graphs with given degrees sequences, \emph{Ann. Comb.}
\textbf{6}, 125-145.\medskip

\noindent Chung, F. and Lu, L. (2002:2), The average distances in
random graphs with given expected degrees, \emph{Proc. Natl. Acad.
Sci.} \textbf{99}, 15879-15882.\medskip

\noindent Dorogovtsev, S. and Mendes, J. (2003), \emph{Evolution
of Networks, from Biological Nets to the Internet and WWW}, Oxford
University Press.\medskip

\noindent Erd\H{o}s, P. and R\'{e}nyi, A. (1959), On random
graphs, \emph{Publ. Math.} \textbf{6}, 290-297.\medskip

\noindent Hofstad, R. van der, Hooghiemstra, G. and Znamenski, D.
(2005), Random graphs with arbitrary i.i.d. degrees, preprint
(www.win.tue.nl/$\sim$rhofstad).\medskip

\noindent Holroyd, A.E. and Peres, Y. (2005), Extra heads and invariant
allocations. {\em The Annals of Probability} {\bf 33}, 31-52.\medskip

\noindent Mattera, M. (2003), Annihilating random walks and perfect
matchings of planar graphs, {\em Discrete Math. and Theor. Comp. Science}
{\bf AC}, 173-180.\medskip

\noindent Meshalkin, L.D. (1959), A case of isomorphisms of Bernoulli schemes,
{\em Dkl. Akad. Nauk. SSSR} {\bf 128}, 41-44.\medskip

\noindent Molloy, M. and Reed, B. (1995), A critical point for
random graphs with a given degree sequence, \emph{Rand. Struct.
Alg.} \textbf{6}, 161-179.\medskip

\noindent  Molloy, M. and Reed, B. (1998), The size of the giant
component of a random graphs with a given degree sequence,
\emph{Comb. Prob. Comp.}\ \textbf{7}, 295-305.\medskip

\noindent Newman, M., Strogatz, S. and Watts, D. (2001), Random
graphs with arbitrary degree distributions and their applications,
\emph{Phys. Rev. E} \textbf{64}, 026118.\medskip

\noindent Parry, W. (1979), An information obstruction to finite expected
coding length, {\em Ergodic Theory} (Proc. Conf., Math. Forschungsinst.,
Oberwolfach 1978), 163-168, Lecture Notes in Math.\ {\bf 729}, Springer,
Berlin.\medskip

\noindent Schmidt, K. (1984), Invariants for finitary isomorphisms
with finite expected coding lengths, {\em Invent. Math.} {\bf 76},
33-40.

\end{document}